\documentclass[9pt]{article}

\usepackage{amssymb}
\usepackage{latexsym}
\input epsf

\def\sref#1{Section~\ref{#1}}

\def\tref#1{Theorem~\ref{#1}}
\def\dref#1{Definition~\ref{#1}}

\def\rref#1{Remark~\ref{#1}}
\def\cref#1{Corollary~\ref{#1}}
\def\pref#1{Proposition~\ref{#1}}
\def\lref#1{Lemma~\ref{#1}}

\def\fgref#1{Figure~\ref{#1}}

\newtheorem{thm}{Theorem}[section] 
\newtheorem{dfn}[thm]{Definition} 
\newtheorem{dfns}[thm]{Definitions} 
\newtheorem{rmk}[thm]{Remark}

\newtheorem{cor}[thm]{Corollary}
\newtheorem{prop}[thm]{Proposition} 
\newtheorem{lem}[thm]{Lemma}
\newtheorem{exs}[thm]{Examples}

\newtheorem{obs}[thm]{Observations}

\newcommand{\intro}[1]
{\renewcommand{\thesection}{\fnsymbol{section}}
\setcounter{section}{-1}
\section{#1}
\renewcommand{\thesection}{\arabic{section}}
}

\def\H#1#2{H^{#1}(#2)}
\def\fr#1#2{\frac{#1}{#2}}

\def\dim#1{{\rm #1}\,}
\def\der#1#2{\frac{\partial #1}{\partial #2}}
\def\sm#1{C^\infty(#1)}

\newcommand{\Pf}{{\em Proof}. }

\newcommand{\cqfd}
{%
\mbox{}%
\nolinebreak%
\hfill%
\rule{2mm}{2mm}%
\medbreak%
\par%
}

\newcommand{\CE}{\mathbb C}
\renewcommand{\d}{{\cal D}{}}  
\newcommand{\e}{{\cal E}{}}
\newcommand{\F}{{\cal F}{}}
\newcommand{\G}{{\cal G}{}}

\newcommand{\p}{{\cal P}{}}

\newcommand{\R}{\mathbb R}

\newcommand{\T}{{\cal T}{}}
\newcommand{\torus}{\mathbb T}

\newcommand{\Z}{\mathbb Z}
\newcommand{\eps}{\varepsilon }
\newcommand{\ol}{\overline}

\newcommand{\co}{cohomology }
\newcommand{\hp}{h-principle }
\newcommand{\tdco}{tangential de Rham cohomology }
\newcommand{\td}{tangential differential }
\newcommand{\s}{symplectic }

\newcommand{\f}{foliation }

\newcommand{\lw}{leafwise }
\newcommand{\al}{almost }
\newcommand{\lsf}{leafwise symplectic form }
\newcommand{\lsfe}{leafwise symplectic form}
\newcommand{\lsst}{leafwise symplectic structure }
\newcommand{\lsste}{leafwise symplectic structure}
\newcommand{\sst}{symplectic structure }
\newcommand{\sste}{symplectic structure}

\newcommand{\lasf}{leafwise almost symplectic form }
\newcommand{\lasfe}{leafwise almost symplectic form}
\newcommand{\lasst}{leafwise almost symplectic structure }
\newcommand{\lasste}{leafwise almost symplectic structure}
\newcommand{\fms}{foliated manifolds }
\newcommand{\fm}{foliated manifold }

\newcommand{\tnd}{tubular neighborhood }

\newcommand{\wrt}{with respect to }
\newcommand{\rp}{respectively }

\begin{document} 

\title{Foliations Associated to Regular Poisson Structures}
\author{M\'elanie Bertelson \thanks{Max-Planck-Institut f\"ur Mathematik, 
D-53111 Bonn, Germany (bertel@mpim-bonn.mpg.de). This work has been supported 
by an Alfred P. Sloan Dissertation Fellowship.}}  
\date{\today}

\maketitle

\begin{abstract}
\noindent 
A regular Poisson manifold can be described as a foliated space carrying 
a tangentially symplectic form. Examples of foliations are produced here 
that are not induced by any Poisson structure although all the basic obstructions 
vanish.       

\end{abstract}

\intro{Introduction}

A {\em Poisson manifold} is a smooth manifold $M$ whose algebra of 
functions $C^{\infty}(M)$ carries a Lie algebra bracket $\{\,,\}$ 
for which each adjoint operator $\{f,\cdot\}$ is a derivation of the 
pointwise multiplication. As revealed by an investigation of its local 
geometry, a Poisson manifold has a natural (possibly singular) foliation 
whose leaves are symplectic manifolds. Conversely, this ``symplectic 
foliation'' determines the Poisson structure. When the foliation is
regular, the symplectic structures on the leaves assemble into a
section of the second exterior power of the cotangent bundle
to the foliation, called a {\em leafwise symplectic structure}.

This work addresses the problem of characterizing regular foliations 
that arise from Poisson structures in this fashion. One motivation for 
this topic is that it constitutes a first step towards classifying
Poisson manifolds (cf.~\cite{AW-Poisson}). More heuristically, any
problem in Poisson geometry potentially splits into two problems~:
one of them concerns foliations; the other one concerns leafwise
symplectic structures. From this perspective, it is important to
understand the relationship between the two. In particular,
investigating the map $\{\mbox{Poisson structures}\} 
\stackrel{\varphi}{\to}\{\mbox{foliations} \}$ is of special relevance. 

The question of existence of symplectic structures on manifolds
appears as the special case where the foliation has a single leaf. For closed
manifolds, this problem is very difficult and mostly unsolved, although
a lot of progress has recently been made. In particular, Taubes, using 
Seiberg--Witten theory, exhibits much subtler obstructions than the basic 
ones (a de~Rham class whose exterior powers do not vanish, and a compatible 
nondegenerate $2$-form)~\cite{CT}. On the other hand, for open manifolds, 
Gromov showed in 1969 that {\em the parametric \hp for open, invariant 
relations is valid} \cite{Gromov1, Gromov-book}. This result implies that 
on an open manifold, every nondegenerate $2$-form is homotopic to a symplectic 
form, and therefore reduces the problem to a question of obstruction theory. 

Our purpose is to understand how far this result can be generalized to
foliated manifolds. In particular, we would like to find a notion of
``openness'' for such spaces that is natural and broad enough, and for
which the following statement holds~: on an open \fm {\em every
leafwise nondegenerate $2$-form is homotopic to a leafwise symplectic
form}. Henceforth, a \fm with noncompact leaves for which this statement
holds will be said to be ``open'' (provided the statement is
nonempty, that is, leafwise nondegenerate $2$-forms exist). The
requirement that the leaves be noncompact is not sufficient to
ensure that the foliation is ``open'', but implies, together with the
existence of a leafwise nondegenerate $2$-form, that each leaf admits
its own symplectic structure. Locally trivial fibrations are ``open'' when 
their fibers are noncompact. This is implied by Gromov's h-principle mentioned 
above.    

Examples of non-``open'' foliated manifolds are produced here. 
Although the foliations exhibited carry leafwise 
nondegenerate $2$-forms, and in some cases display no identifiable
closedness feature, they do not carry {\em any} leafwise symplectic
form. An intriguing example consists of a foliation associated to a
submersion whose fibers are all noncompact, connected, and
diffeomorphic to one  another. It is not a locally trivial fibration
of course, but seems very close to being one. These examples show that
one cannot be too optimistic about the size of the class of ``open''
foliated manifolds, and they illustrate different types of simple
obstructions that one may encounter when trying to build leafwise
symplectic structures. 

In a subsequent paper \cite{B-hp} (see also~\cite{B}), the following
theorem is proved~: {\em On a foliated manifold $(M,\F)$ that admits 
a positive, proper function $f:M\to\R$, without leafwise local maxima, 
satisfying certain generic nondegeneracy conditions, the \hp for open 
and foliated-invariant relations is valid.} In particular, under the
hypotheses of this theorem, every leafwise nondegenerate $2$-form is
homotopic to a leafwise symplectic form. The condition of existence of
such a function seems thus to be a good candidate for the definition
of {\em open foliated manifold}, all the more, since it corresponds
to the usual definition in the case of foliations with a single leaf.
Besides, it is not difficult to verify directly that the foliated 
manifolds described here do not support proper functions without 
leafwise local maxima.\\

{\bf Acknowledgements.} I wish to thank Alan Weinstein for suggesting 
this problem and for constant and helpful advice, as well as Yasha 
Eliashberg for many insightful conversations.

\section{Preliminaries}\label{blabla}

\subsection{Regular Poisson manifolds}

The definition of regular Poisson manifolds, and the existence of the 
associated symplectic foliation are recalled in this section.

\begin{dfn} A {\em Poisson manifold} is a smooth manifold whose 
algebra of functions $C^\infty(M)$ is endowed with a Lie 
bracket $\{ \,,\}$ acting on each argument as a derivation, that 
is, the following identity (called Leibniz identity) is 
satisfied for all $f$, $g$ and $h$ in $C^\infty(M)$~:
$$
\{fg,h\}=g\{f,h\}+f\{g,h\}. 
$$
\end{dfn}
As a consequence of this definition, a vector field $X_f = \{f,\cdot\}$ 
is associated to each function $f \in C^\infty(M)$. Moreover, the map 
$C^\infty(M) \to \Gamma TM : f \mapsto X_f$ is a Lie algebra homomorphism.  

\begin{dfn}
When the rank of the distribution $\d = \{X_f(x) ; f\in C^\infty(M), 
x\in M\}$ is constant, the Poisson manifold is said to be {\em regular}. 
\end{dfn}
A regular Poisson manifold supports a foliation whose leaves are symplectic 
manifolds. Indeed, the distribution $\d$ is involutive, thus integrable, and 
a leaf $F$ carries the symplectic structure $\omega_F$ defined by the 
following expression~:
$$
\omega_F (X_f\rule[-2mm]{0.1mm}{5mm}_F,X_g\rule[-2mm]{0.1mm}{5mm}_F) = 
\{df,dg\}\rule[-2mm]{0.1mm}{5mm}_F \quad f, g \in C^\infty(M).
$$
The symplectic structures along the leaves assemble into a global 
section of the bundle $\Lambda^2 T^*\F$. (Here $T\F$ denotes the tangent 
bundle to the foliation, and $T^*\F$ denotes its dual.) Conversely, 
the data of a foliation $\F$, endowed with a section $\omega$ of the 
bundle $\Lambda^2 T^*\F$ whose restriction to each leaf is a symplectic 
form, determines a regular Poisson structure having $\F$ as associated 
foliation.  

\subsection{Tangential de Rham cohomology}

This section introduces natural extensions to foliated spaces of the 
complex of differential forms and de Rham cohomology. This material can 
be found in \cite{HMS} for example. \\

Let $\Omega^k(\F)$ denote the space of sections of the bundle 
$\Lambda^k T^* \F$. The elements of $\Omega^k(\F)$ are called {\em tangential 
differential $k$-forms}, and can be thought of as skewsymmetric 
$\sm{M}$-multilinear maps $\Gamma T\F \times \ldots \times \Gamma T\F 
\to  C^\infty (M)$. As for ordinary differential forms, the expression   
\begin{eqnarray*}	
d_\F \alpha (X_0,\ldots,X_k) & = & \sum^k_{i=0} (-1)^i X_i \left(\alpha
(X_0,\ldots,\hat X_i,\ldots,X_k)\right) + \\ &&\sum_{i<j}(-1)^{i+j} 
\alpha([X_i,X_j], X_0,\ldots,\hat X_i,\ldots,\hat X_j,\ldots,X_k) \;,
\end{eqnarray*}
where $\alpha\in \Omega^k(\F)$, where $X_0,\ldots,X_k \in \Gamma T\F$, and
where a symbol covered by a hat is omitted, defines a degree one differential 
operator $d_\F$ that satisfies $d_\F^2 = 0$. A tangential differential 
$k$-form $\alpha$ with $d_\F\alpha=0$ is said to be {\em $d_\F$-closed}, 
and if $\alpha=d_\F\beta$ for some tangential differential $(k-1)$-form 
$\beta$, then $\alpha$ is said to be {\em $d_\F$-exact}. The induced \co 
$$
H^\star(\F)=\frac{\{\mbox{$d_\F$-closed $\star$-forms}\}}{ \{\mbox{ 
$d_\F$-exact $\star$-forms}\}}
$$ 
is called the {\em tangential de~Rham cohomology}. \\  

The natural inclusion $i:T\F\to TM$ induces surjective bundle maps $r_k : 
\Lambda^k T^*M \to  \Lambda^kT^*\F$. In particular, any tangential differential 
form $\alpha$ is the restriction to $\Gamma T\F \times \ldots \times \Gamma 
T\F$ of some ordinary differential form, called hereafter an {\em extension 
of $\alpha$}. 

\begin{rmk}{\rm If the foliation has a single leaf, we recover
the ordinary de~Rham cohomology.
}\end{rmk} 
  
The following two properties (\pref{M-V} and \pref{Kunneth-prop}) will 
be useful in the next section, where calculation of the \tdco for some 
specific foliations will be needed.
 
\begin{prop}[\cite{ElKacimi}]\label{M-V}(Mayer--Vietoris sequence.) Let 
$(M,\F)$ be a foliated manifold, and let $U$ and $V$ be open subsets of
$M$. Then the following short sequence is exact~:
$$\begin{array}{ccccccccc}
0 & \to & \Omega^k(\F|_{U\cup V}) & \to & \Omega^k(\F|_{U}) \oplus 
\Omega^k(\F|_{V}) & \to & \Omega^k (F|_{U\cap V}) &\to & 0\\
& & \alpha & \mapsto & (\alpha|_U,\alpha|_V)&&&&\\
& & & & (\beta,\gamma) & \mapsto & \gamma|_{U\cap V} - \beta |_{U\cap 
V} & \;.& 
\end{array}$$
It therefore induces a long exact sequence in cohomology~:
\small{
$$
\ldots \to H^k(\F|_{U\cup V}) \to H^k(\F|_{U})\oplus 
H^k(\F|_{V}) \to H^k (F|_{U\cap V}) \to H^{k+1}(\F|_{U\cup 
V}) \to \ldots
$$
}
\end{prop}
The proof is similar to that of exactness of the ordinary Mayer--Vietoris 
sequence in de~Rham cohomology.\\

Let $(M_1,\F_1)$ and $(M_2,\F_2)$ be two foliated manifolds, and consider 
their Cartesian product $(M_1\times M_2, \F_1\times\F_2)$. For every pair 
$(i,j)$ of integers, there is a pairing  

$$\begin{array}{ccccl}
\times_{i\,j} & : & \Omega^i(\F_1)\otimes \Omega^j(\F_2) & \to & 
\Omega^{i+j}(\F_1 \times \F_2)\\
& & (\alpha_1 \otimes  \alpha_2) & \mapsto & p_1^* (\alpha_1) \wedge
p_2^* (\alpha_2)\;,
\end{array}$$ 
where $p_\ell$ denotes the natural projection of $M_1\times M_2$ onto its
$\ell$-th factor. These pairings induce maps in cohomology~:   
$$
\times^k : \bigoplus_{i+j=k}H^i(\F_1)\otimes H^j(\F_2) \to H^{k}(\F_1
\times \F_2)\;.
$$
When the map $\times^k$ is an isomorphism, we say that the K\"unneth 
formula  
$$
H^{\star}(\F_1 \times \F_2)  \simeq \bigoplus_{i+j=\star} H^i(\F_1) \otimes
H^j(\F_2)
$$ 
is valid. The following result will be sufficient for our needs.  

\begin{prop}[\cite{B}]\label{Kunneth-prop}
Let $(M_1,\F_1)$ be any foliated manifold, and let $M_2$ be a manifold 
of finite type foliated by a single leaf. Then the K\"unneth formula
associated to the product $(M_1,\F_1) \times M_2$ is valid, that is,
$$
H^{k}(\F_1 \times M_2) \simeq \bigoplus_{i+j=k} H^i(\F_1) \otimes
H^j(M_2) \quad \mbox{for all $k$.}
$$
\end{prop}
The proof is again very similar to the proof of the corresponding result 
in the single-leaf case.   
 
\subsection{Terminology and first examples}

\begin{dfns}\label{ls} Let $(M,\F)$ be a foliated manifold.
\begin{enumerate}
\item[-] A {\em \lasf} is a tangential differential $2$-form $\beta$ that 
is nondegenerate, that is, the bundle map $T\F \to T^*\F : X \mapsto 
\beta(X,\cdot)$ is an isomorphism.  
\item[-] A {\em \lsf} is a $d_\F$-closed, nondegenerate, \td $2$-form. 
Equivalently, it is a tangential differential $2$-form whose restriction 
to each leaf is a symplectic form.  
\item[-] A {\em \lw volume form} is a nonvanishing tangential differential 
form whose degree coincides with the dimension of the foliation. The \fm 
$(M,\F)$ is said to be {\em orientable} if it admits a \lw volume form.   
\end{enumerate}
\end{dfns}
With this terminology, the data of a regular Poisson structure on a manifold 
is equivalent to that of a foliation endowed with a \lsfe.
 
\begin{exs}{\rm (Examples of regular Poisson structures.)
\begin{enumerate}
\item[a)] A symplectic manifold is a regular Poisson manifold. Any
          manifold endowed with the trivial Poisson bracket $P=0$ is
          Poisson manifold as well. 
\item[b)] Let $M$ be a manifold endowed with a two-dimensional foliation 
          $\F$, and suppose that $(M,\F)$ is orientable. Then any \lw 
          volume form is a leafwise symplectic form (and vice
          versa). The set of \lsste s is therefore parameterized by the
          set of nonvanishing smooth functions.
\item[c)] Let $(M,P)$ be any Poisson manifold. The rank of $P$ is 
          maximum on $O_1$, an open subset of $M$. Similarly, the 
          rank of $P|_{M-\ol O_1}$ is maximum on an open subset $O_2$,
          and so forth. This shows that $M$ contains a dense open set 
          $O=\cup_iO_i$ such that $(O,P|_O)$ is a disjoint union of 
          finitely many regular Poisson manifolds.
\item[d)] As explained in the next section, a locally trivial fibration 
          with noncompact fibers that admits a \lasf also admits a \lsfe.
\item[e)] More examples can be constructed by means of the criterion for 
          existence of \lsst exhibited in \cite{B-hp}.
\end{enumerate}      
}\end{exs} 

\begin{rmk}\label{Anosov}{\rm 
As illustrated by the following example, a compact manifold may support 
$d_\F$-exact \lsfe s. Let $A$ be an element in $SL(2,\Z)$ with ${\rm 
trace}\,(A)>2$. Consider the \fm $(M,\F)$ constructed by suspension of $A$, 
thought of as an isomorphism of the \fm $(\torus^2, \F_X)$, where $\F_X$ 
denotes the (irrational slope) linear foliation of $\torus^2$ in the 
direction of an eigenvector $X$ of $A$. It is shown in \cite{ElKacimi} 
that the second \tdco space of $\F_X$ vanishes; the \lsste s on $(M,\F)$ 
are therefore all exact.   
}\end{rmk} 

\subsection{Description of the problem}
 
Given a manifold $M$, consider the map 
$$\begin{array}{c}
\{ \mbox{ rank-$2k$ regular Poisson structures on $M$} \} \\
\downarrow \varphi \\
\{ \mbox{ $2k$-dimensional foliations on $M$} \}\;.
\end{array}$$
As explained in the introduction, the problem approached in this work 
is that of describing the image of this map. More specifically, we aim 
at finding interesting necessary and sufficient conditions on a
foliation $\F$ to be in the image of $\varphi$. Along those lines, 
the following observations are relevant.  

\begin{obs}\label{Tintoretto}{\rm \hspace{-.1cm}

(I) A foliation $\F$ supporting a \lw \s form obviously satisfies the
following conditions.
\begin{enumerate}
\item[i)] Each leaf of $\F$ supports some symplectic form.
\item[ii)] The \fm $(M,\F)$ supports some \lw almost symplectic form
(\dref{ls}). Equivalently, $T\F$ admits a structure of symplectic
vector bundle. In particular, the \fm $(M,\F)$ is orientable. 
\end{enumerate}

(II) For foliations with a single leaf, the problem becomes that of
characterizing manifolds supporting some symplectic form. There is a
strong dichotomy between open and closed manifolds. For closed 
manifolds, the existence problem is hard and essentially unsolved. 
The basic obstructions are~:  
\begin{enumerate}
\item[a)] Existence of a class $a$ in the second de Rham \co space of $M$ such 
          that $a^k\neq 0$, where $k = \fr 1 2 \,\dim{dim} M$.
\item[b)] Existence of an \al \s structure or equivalently of an almost complex 
          structure inducing a volume compatible with $a^k$.   
\end{enumerate}
 
Results of Taubes involving Seiberg--Witten theory imply that
fulfillment of these two conditions does not guarantee existence of a
\s form. For instance, the connected sum of three copies of $\CE P^2$
satisfies conditions a) and b), but does not carry any symplectic form
(\cite{CT}). For open manifolds, when one does not require the
symplectic structure to have any particular behavior at infinity or at
the boundary, the problem has been solved by Gromov in his thesis. He
proved the following result.       

\begin{thm}[\cite{Gromov1}]\label{Gromov}
On an open manifold $M$, any almost symplectic form is homotopic
among almost symplectic forms to a symplectic form. 
\end{thm}

(III) The next simplest foliations are locally trivial fibrations. 
\tref{Gromov} (or rather a stronger version of this theorem) implies the 
following result.  

\begin{cor}[\cite{B}]\label{fib}
Let $\pi : M \to B$ be a locally trivial fibration whose fiber is an
open manifold $F$. Then any \lw almost symplectic form on $(M,\F)$ is
homotopic, among such forms, to a \lw symplectic form.  
\end{cor}

(IV) In view of \tref{Gromov}, if the foliation at hand 
has no compact leaf, then property ii) above implies property i). 
}\end{obs}

Our approach has been to focus on foliations with noncompact leaves,
and search for a generalization of \tref{Gromov} alongside guiding
examples. It is important to keep in mind that one does not control
the behavior at infinity of the \sst constructed in the proof of
\tref{Gromov}. In contrast, a symplectic structure on a leaf of a
foliation that extends to a \lsst on the entire \fm is most likely very
constrained at infinity, due partly to recurrence phenomena, partly to
the influence of neighboring symplectic leaves. We have been trying to
grasp the type of behavior at infinity that a foliation supporting a
\lsst typically has. The examples in \sref{examples} provide
some answers to natural questions directed along those lines. 
A generalization of \tref{Gromov} will appear in \cite{B-hp}. 
Its statement is as follows~:    

\begin{thm}[\cite{B-hp}]\label{pre-Stravinsky} Let $(M,\F)$ be an {\em open
foliated manifold}. Then any \lasf is homotopic, among \lasfe s, to a
leafwise symplectic form.
\end{thm}
An {\em open \f manifold} is defined to be \fm that admits a positive, proper 
function, without leafwise local maxima, whose jet satisfies some 
transversality condition (see \cite{B-hp, B} for a precise definition).  

\begin{rmk}\label{Braque}{\rm Let $(M,\F)$ be a foliated manifold that 
admits a proper function $f : M\to [0,\infty)$ without leafwise local maxima. 
Consider a partition $\p = \{a_0 = 0 < a_1 < \ldots < a_k < \ldots \}$ 
of $[0, \infty)$ by noncritical values of $f$, and the associated exhaustion 
of $M$ by the slices $K_i = f^{-1}([a_i,a_{i+1}])$. One observes that the 
condition that $f$ has no \lw local maxima implies that for any leaf $F$ 
of the foliation $\F$, and for any index $i$, the set $F -\dim{int}(K_i)$ 
has no compact connected component ($f|_F$ would necessarily achieve a
local maximum on such a component).   
}\end{rmk}
Using this criterion, it is easy to verify directly (that is, independently 
of \tref{pre-Stravinsky}) that the foliated manifolds constructed in the next 
section are not open.   

\section{Examples}\label{examples}
 
This section could be entitled {\em variations on the
theme :}
\begin{center}\label{Giotto}
{\em a compact manifold with vanishing second de~Rham \co} \\ 
{\em is never symplectic.} 
\end{center}
It presents examples of foliated manifolds with noncompact 
leaves that do not carry {\em  any} \lsste, although \lasste s 
exist. For all of them, the second \tdco space
vanishes, so that a \lsst would necessarily be exact, and the argument
consists in proving that the \fm at hand may not carry any exact
\lsste. In each case, the leaves being noncompact, the argument is not as
simple as that of the {\it theme}, but presents nevertheless a strong
analogy with it. Therefore, it seems tempting, in view of \tref{Gromov},
to attribute nonexistence of exact \lsste s to some kind of compactness,
or non-openness, in disguise. We have followed this line of ideas
through the exposition of the examples, trying to provide some
interpretation for each. 

The exposition is organized as follows. The examples are presented
along with arguments for nonexistence of exact \lsste s. There is not a
one to one correspondence between the two. The arguments are
independent, but the examples are not really so. They evolve from
seemingly quite closed to seemingly more and more open, so that, as
already mentioned earlier, the last one is the most seemingly open one
and therefore the most interesting one too. Nevertheless, we felt it
worthwhile to include the previous examples as~well.

\subsection{Leafwise almost symplectic structures}\label{laze}

We begin with showing a simple conditions ensuring existence of a 
\lasst (\dref{ls}). 

\begin{lem}\label{cond1} A \fm of the type $(N,\G) \times P$, where $\G$
is the foliation of $N$ by the orbits of a nonvanishing vector field
$X$, and where $P$ is an odd-dimensional manifold that admits a
contact form, always supports some \lasste.  
\end{lem}
Recall that a {\em contact form} is a $1$-form $\alpha$ such that $\alpha
\wedge (d \alpha)^k$ does not vanish for $k = \fr 1 2 (\dim{dim} P
- 1)$.\\

\Pf Let $\beta$ be a $1$-form on $N$ such that $\beta(X)=1$. Then 
$\beta \wedge \alpha + d \alpha$ extends a \lw \al \s form.
\cqfd

More generally, define a {\em \lw \al contact structure} on a 
manifold  endowed with a $(2k+1)$-dimensional foliation to be a pair
$(\alpha,\beta)$, consisting of a tangential differential $1$-form
$\alpha$ and a tangential differential $2$-form $\beta$, such that
$\alpha \wedge (\beta)^k$ is \lw volume form. If $(M_1,\F_1)$  and
$(M_2,\F_2)$ admit \lw \al contact structures $(\alpha_1,\beta_1)$ and
$(\alpha_2,\beta_2)$ respectively, then $\beta_1 + \alpha_1 \wedge
\alpha_2 + \beta_2$ is a \lasst on $(M_1 \times M_2 , \F_1 \times
\F_2)$.  

\subsection{First variation}

Let $a$ be an irrational number, and consider the \fm $(\torus^2,
\F_a) \times S^3$, where $\F_a$ denotes the linear foliation of $\torus^2$ 
of slope $a$. It admits a \lasst (\lref{cond1}), but not any \lsste. 
The latter statement is, as will be seen below, a consequence of the 
following proposition. 

\begin{prop}\label{var1}
A \fm $(M,\F)$ satisfying the following two conditions does not 
support any exact \lsste. 
\begin{enumerate}
\item[-] $M$ is compact,
\item[-] $(M,\F)$ admits a transverse volume form $\mu$.
\end{enumerate}
\end{prop}
Let us recall that a transverse volume form on a manifold $M$, endowed
with a foliation $\F$ of 
codimension $q$, is a nonvanishing closed $q$-form $\mu$ such that 
$i(X)\mu=0$ for any vector field $X$ tangent to $\F$. It is easily seen 
that the \fm $(\torus^2, \F_a) \times S^3$ admits such a form. 
(If $\theta$ and $\varphi$ are coordinates on $\torus^2$ for which 
the vector field $\der{}{\theta} + a \der{}{\varphi}$ is tangent to 
the foliation $\F_a$, the form $$\mu = p_1^*(a\, d\theta - d\varphi)\;,$$ 
where $p_1$ is the natural projection of $\torus^2 \times S^3$ onto its first
factor, is a transverse volume form.) Observe that, as implied by 
\pref{Kunneth-prop}, $\H2{\F_a \times S^3} = 0$. (Notice that the manifold 
$S^3$ could be replaced by any homology three-sphere, or any odd-dimensional 
manifold $P$ endowed with an almost contact structure and satisfying $H^1(P) 
= 0 = H^2(P)$).\\ 

\pagebreak

{\em Proof of \pref{var1}.} Suppose on the contrary that $(M,\F)$
admits an exact \lsst $\omega = d_\F \alpha$. Let $\tilde\omega \in
\Omega^2(M)$ be an extension of $\omega$. Then the form
$(\tilde\omega)^k \wedge \mu$ is a volume form on $M$. We claim that
$(\tilde\omega)^k \wedge \mu$ is exact. Indeed, $(\tilde\omega)^k
\wedge \mu = d (\tilde\alpha \wedge (\tilde \omega)^{k-1} \wedge
\mu)$, where $\tilde\alpha$ is an extension of~$\alpha$. 
\cqfd

\subsection{Second variation}

In the previous example, closedness of the ambient manifold plays, it
seems, a crucial role. Besides, it is to be expected that a foliation
on a closed manifold presents some features of a compact nature, and
that a definition of openness for foliations would require the ambient 
manifold to be nonclosed. Nevertheless, compactness of the ambient 
manifold is not, as will be seen, the only missing obstruction to 
existence of a \lsste. \\ 

Starting from the preceding example $(\torus^2,\F_a) \times S^3$
and removing a point $p$ from $\torus^2$ provides a noncompact \fm
$(M,\F) = (\torus^2 - \{p\},\F_a|_{\torus^2 - \{p\}}) \linebreak
\times S^3$ 
with open leaves (\fgref{fig:punct-torus}) that admits a \lasst 
(since $(\torus^2,\F_a) \times S^3$ already did). Nevertheless, 
the \fm $(M,\F)$ does not support any \lsste. This is a consequence 
of the following proposition, combined with the fact that $H^2(\F)$ 
vanishes (as implied by \pref{Kunneth-prop}).  

\begin{figure}[h]
\begin{center}
\input punct-torus.pstex_t
\smallskip\noindent 
\caption{The foliated manifold $(\torus^2 - \{p\},\F_a|_{\torus^2 -
\{p\}}) \times S^3$} 
\label{fig:punct-torus}
\end{center}
\end{figure}

\begin{prop}\label{var2}A \fm $(M,\F)$ of the type $(N,\G) \times P$,
where $P$ is any closed odd-dimensional manifold, and where $\G$ is the
foliation of $N$ by the orbits of a nonvanishing vector field $X$
admitting some nonclosed orbit, may not support any exact \lsste. 
\end{prop}

\Pf Suppose that $\alpha$ is a \td $1$-form on $(M,\F)$ for which
$d_\F \alpha = \omega$ is a \lsste. Let $2k =\dim{dim}\F$. If $i : S
\to M$ is an immersed oriented submanifold of dimension $2k-1$ contained
in a leaf $F$, define   
$$
\nu (S) = \int_S i^*(\alpha \wedge \omega^{k-1})\;.
$$
Observe that if $S$ bounds a compact domain $D$ in $F$, then 
$$\nu(S) = \pm \int_D \omega^k\;.$$

We show that the  symplectic volume of a nonclosed leaf $F$ of $\F$ 
is necessarily infinite. 
Let $a$ be an accumulation point of the leaf $F$ that does not lie in
$F$, and let $(U,\varphi)$ be a chart adapted to $\F$ and centered at
$a$. For $x$ in $U$, denote by $P_x$ the plaque of $\F|_U$ containing
$x$. Consider a sequence $(a_i)$ contained in $U\cap F$ that converges
to $a$. Let $\delta$ be the \s volume of the plaque $P_a$. Then, for
$i$ sufficiently large, the \s volume of the plaque $P_{a_i}$ is at
least~$\fr{\delta}2$. We have thus infinitely many disjoint subsets of
$F$, each which contributes at least~$\fr{\delta}2$ to the volume
of~$F$. The latter must therefore be infinite. Notice that this 
argument shows that the volume of the end $\e$ of $F$ containing
the sequence $(a_i)$ is infinite as well. \\

On the other hand, the volume of $\e$, or rather of a set $U \simeq 
[0,\infty) \times P$ representing $\e$, can be written as follows~: 
\begin{eqnarray*}
\lim_{n\to\infty} \int_{[0,n]\times P} \fr{\omega^k}{k!} & = &
\lim_{n\to\infty} \int_{\partial [0,n]\times P} \fr1{k!} \alpha \wedge
\omega^{k-1} \\ 
& = & \lim_{n\to\infty} \fr1{k!}\left(\nu(\{n\} \times P) - \nu(\{0\}
\times P)\right) \;,
\end{eqnarray*}
where $\{t\} \times P$ is oriented as the boundary of $[0,t] \times P$, 
the latter being oriented by the \lw symplectic form $\omega$. Because
the volume of $\e$ is infinite, $\lim_{n\to\infty}\nu(\{n\} \times P) = 
\infty$. But the map   
$$
N \to \R : q \mapsto \nu(\{q\} \times P)
$$
is continuous (even smooth). Therefore, if $q_a$ denotes the first
coordinate of $a$ \wrt the splitting $M = N \times P$, the quantity
$\nu(\{q_a\} \times P)$ has to be infinite, and we reach a contradiction. 
\cqfd

\begin{rmk}{\rm Observe that if some orbit $\varphi^\R(x)$ of the vector 
field $X$ is contained in its positive (or negative) limit set 
$\lim_{t\to \infty} \varphi^t(x)$ (as is the case for a vector field 
tangent to the foliation $(\torus^2 - \{p\},\F_a|_{\torus^2 - \{p\}}$)), 
the conclusion of \pref{var2} can be strengthened as follows~: the foliation 
$(M,\F)$ may not support any \lsst whose restriction to the leaf $F = 
\varphi^\R(x) \times P$ is an exact form. Therefore, to obtain 
nonexistence of \lsste s on such a foliation, it is sufficient that the 
second de~Rham \co of the leaf $F$ vanishes. In particular, the \fm 
$(\torus^2 - \{p\}, \F_a|_{\torus^2 - \{p\}}) \times P$, where $P$ is a 
manifold with $H^2(P) = 0$ does not carry any \lsste. 
}\end{rmk}

\begin{rmk}{\rm The above lemma also proves that the first example 
$(\torus^2,\F_a)\times S^3$ does not carry any \lsste. More generally, 
it shows that a \fm of the type considered in the previous proposition 
may not support any \lsst whose top exterior power is exact (simply replace 
$\alpha \wedge \omega^{k-1}$ in the proof of \pref{var2} by a primitive 
of $\wedge^k \omega$). 
This argument may
also be generalized in another direction, namely to \fms having a leaf
with a periodic end, that is, an end diffeomorphic to an
infinite connected sum $P\#P\#\ldots$, where $P$, the {\em period}, is a
manifold whose boundary splits into two diffeomorphic pieces. One has to
pay attention, though, to the fact that the hypersurfaces along which
the connected sum is performed (the $\{n\} \times P$'s in the previous
proof) might not ``converge'' in $M$ to a hypersurface in a
leaf. Indeed, looking back at the Anosov foliation, described in
\rref{Anosov}, one notices that its leaves are either planes or
cylinders. Thus they have nonclosed periodic ends. Nevertheless, this
foliation admits \lsste s although its second \tdco space vanishes. The
argument used in the proof of \pref{var2} must therefore fail. 
}\end{rmk}

In the previous example, the leaves were not embedded, and it is understandable 
that this might create obstructions to existence of \lsste s. Indeed, a 
nonembedded leaf carries an additional topology $\T$, that is coarser than 
the topology associated to its manifold structure. A \s form on such a 
leaf that is the restriction of a \lw symplectic form would necessarily 
be in some sense continuous \wrt $\T$. It could just happen that some of 
the nonembedded leaves do not admit any \s structure possessing this additional 
property (this is the case for a leaf of the foliation $(\torus^2 - 
\{p\},\F_a|_{\torus^2 - \{p\}}) \times S^3$). Nevertheless, \pref{var2} does 
not rely on the leaves being nonembedded, it also applies to some foliations 
with nonclosed embedded leaves, as illustrated by the following example. 
Thus, the presence of an additional topology is not the only problem either.

\begin{figure}[h]
\begin{center}
\input embedded-punct-torus.pstex_t
\smallskip\noindent 
\caption{The foliated manifold $(\torus^2 - \{p\}, \F) \times S^3$}
\label{fig:embedded-punct-torus}
\end{center}
\end{figure}

Consider the \fm $(M,\F) = (\torus^2 - \{p\}, \F) \times S^3$, where
the foliation $\F$ is described as follows. Let $\theta$ and $\varphi$
be standard coordinates on $\torus^2$ vanishing at $p$. Consider the
vector field $X = \der{}{\varphi} + f \der{}{\theta}$ where $f$ is a
smooth function on $\torus^2$ such that  
\pagebreak
\begin{enumerate}
\item[-] $f$ depends on $\theta$ only,
\item[-] $f(0)=0$,
\item[-] $f(\theta) > 0$ for $\theta \neq 0$.
\end{enumerate}
Let $\G$ denote the foliation of $\torus^2$ by the orbits of $X$. It has
one closed leaf ($\theta = 0$). The nonclosed leaves are embedded and
accumulate on that closed leaf. The foliation we are interested in is
$\F = p_1^*(\G|_{\torus^2-\{p\}})$ (see \fgref{fig:embedded-punct-torus}), 
where $p_1$ denotes the projection $(\torus^2 - \{p\}) \times S^3 \to \torus^2 - 
\{p\}$. Because the point $p$ has been removed, the leaves
of the foliation $\F$ are all noncompact. Here again, $H^2(\F)$
vanishes (\pref{Kunneth-prop}), and $(M,\F)$ supports some \lasst
(\lref{cond1}). However, \pref{var2} implies that $(M,\F)$ does not
support any \lsste.   

\subsection{Third variation and main example}\label{Berlioz}

We will present here an example of a foliation given by a 
submersion whose fibers are open, connected, and diffeomorphic to one 
another, that admits \lasste s but not any \lsste. This foliation is 
not, of course, a locally trivial fibration (cf.~\cref{fib}), but 
seems nevertheless very close to being one. \\

\begin{figure}[h]
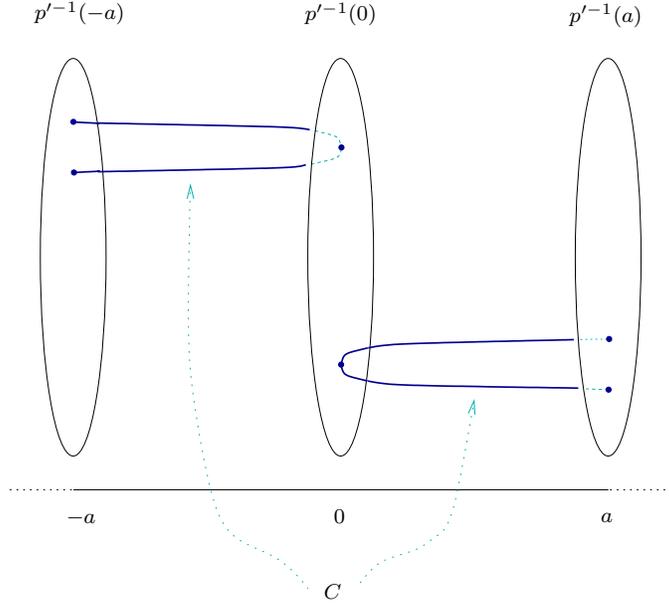

\begin{center}
\input Brahms.pstex_t
\smallskip\noindent 
\caption{The curve $C$ sitting in $M'$}
\label{fig:Brahms}
\end{center}
\end{figure}

Let $F$ be any closed manifold of dimension at least four with $H^2(F)=0$
that supports an almost symplectic structure (take $F = S^1\times S^3$
for instance). Consider the manifold $M' = F\times (-1,1)$, the
projection $p': M' \to (-1,1)$, and the foliation $\F'$ by the fibers of
$p'$. We will remove from $M'$ a closed set $C$, consisting of two
disjoint embedded lines, that intersects each leaf of $\F'$ along
exactly two points. Explicitly, let $p_1, q_1, p_2, q_2$ be four
distinct points in $F$. Suppose that $p_1$ and $q_1$ (\rp $p_2$ and 
$q_2$) are contained in the interior of an embedded closed ball $U_1$ 
(\rp $U_2$), and that $U_1 \cap U_2 = \emptyset$. For $i=1,2$, let 
$\gamma_i : (-1,1) \to U_i$ be a path joining $p_i$ to $q_i$. Define $C$ 
as follows~:  
$$
C = \left(\bigcup_{t\in (-1,0]}\{\gamma_1 (t),\gamma_1(-t)\} \times \{t\}
\right) \;\; \cup \;\; \left(\bigcup_{t\in [0,1)} \{\gamma_2(t), 
\gamma_2(-t)\}\times\{t\}\right)\;.
$$
The foliated manifold $(M=M'- C, \F= \F'|_{M' - C})$ (see
\fgref{fig:Brahms}) admits a foliated almost symplectic structure. Its
leaves are open, hence symplectic. But it does not support any \lsste. The
proof below of the last assertion relies mostly on the vanishing of
the second tangential de~Rham cohomology space of $\F$.  

\begin{lem}\label{Botticelli}
$$H^2(\F) = 0\;.$$
\end{lem}

\Pf For $i=1,2$, let $\tau_i$ be an open \tnd of $\partial U_i$ in $F$
such that if $U'_i= U_i \cup \tau_i$ then $U'_1\cap U'_2 = \emptyset$. 
Consider $W = M \cap \left((U'_1\cup U'_2) \times (-1,1) \right)$ and
$W'= M - ((U_1 \cup U_2)\times (-1,1))$. The Mayer--Vietoris
sequence for the open subsets $W$ and $W'$ (\pref{M-V}) contains this part 
$$
\ldots \to H^1(\F|_{W\cap W'}) \to H^2(\F) \to 
H^2(\F|_{W})\oplus H^2(\F|_{W'}) \to H^2(\F|_{W\cap W'}) \to \ldots
$$
Since $\F|_{W\cap W'}$ and $\F|_{W'}$ are trivial foliations with 
leaves diffeomorphic to $(S^3 \coprod S^3) \times \R$ and $F - (U_1\cup 
U_2)$ respectively, $H^i(\F|_{W\cap W'})={0}$ for $i=1,2$, and  
$H^2(\F|_{W'})={0}$ (cf.~Lemma in \cite{GLSW}). Hence $H^2(\F)$ is  
isomorphic to $H^2(\F|_W)$. On the other hand, the \fm $(W, \F|_W)$ 
is made of two disjoint pieces isomorphic to the following model. \\ 

Let $B^n$ (\rp $B^n_r$) denote the open ball in $\R^n$ of radius $1$ 
(\rp $r$). For $n\geq 4$, consider $N=B^n \times (-1,1) - C$, where 
$C = \cup_{t \in [0,1)} \{\gamma(t), \gamma(-t)\} \times \{t\}$, and 
where $\gamma : (-1,1) \to B^n$ is a smooth path joining two distinct 
points in $B^n$. Suppose for convenience that $\gamma(0)=0$. Endow $N$ 
with the foliation $\G$ by the fibers of the natural projection $\pi : 
N \to (-1,1)$.  

\begin{lem}\label{quatuor} $$H^2(\G)=0\;.$$
\end{lem}

\Pf Let $\beta$ be a $d_\G$-closed \td $2$-form on $(N,\G)$. Since
$\G$ is trivial on $U^-=\pi^{-1}((-1,0))$ and on $U^+=\pi^{-1}((0,1))$,
the Lemma in \cite{GLSW} implies that there exists \td $1$-forms 
$\alpha^-$ and $\alpha^+$, defined on $U^-$ and $U^+$ respectively, 
such that
$$
d_\G \alpha^\pm = \beta|_{U^\pm}\;.
$$
Let $\phi_1$ and $\phi_2$ be smooth functions on $N$ such that
$$
\begin{array}{cc}
\phi_i = & \left\{ \begin{array}{lll}
             1 & \mbox{near} & \overline{O_{2i}} \\
             0 & \mbox{near} & N - O_{2i-1}\;, 
             \end{array} \right. 
\end{array}
$$
where $O_1, O_2, O_3, O_4$ are open subsets of $N$ defined as follows~: 
$$
O_j = N \cap \left(\bigcup_{t\in(-1,1)} B^n_{\frac{|t|}{j}}(0)\times
\{t\}\right)\;. 
$$
We suppose that 
$$
\{\gamma (t), \gamma (-t)\} \subset B^n_{\frac{|t|}{4}}(0)\;\; \mbox{
for all } \;\; t\in (-1,1)\;.
$$ 
(If this is not true, replace in the definition of the $O_j$'s the
function $t \mapsto |t|$ by an appropriate continuous function vanishing 
at $0$). Define    
$$
\alpha = \left\{\begin{array}{ccl} 
           \phi_1 \alpha^- & \mbox{on} & U^- \\
           \phi_1 \alpha^+ & \mbox{on} & U^+ \\
           0 & \mbox{on} & \pi^{-1}(0)\;.
           \end{array}\right.
$$
By construction, $\alpha$ is a smooth section of $T^*\G$, and the form 
$\beta' = \beta - d_{\G}\alpha$ vanishes near $\ol{O_2}$. Let $\alpha'$ 
be a \td $1$-form on $U = N - \overline{O_4}$ such that $d_\G \alpha' 
= \beta'|_{U}$. The form $\alpha'$ exists because the foliation $\G|_U$ 
is isomorphic to a product foliation whose leaves are diffeomorphic to 
$B^n-\{0\}$. Since $\alpha'$ is $d_\G$-closed on the open set 
$\tilde{U}= (O_2 - \overline{O_4})$, we have $\alpha'|_{\tilde{U}} =
d_\G f$ for some smooth function $f$ defined on $\tilde{U}$. Let
$$
\alpha'' =\left\{ \begin{array}{cl}
            \alpha' \, -\, d_\G(\phi_2 f) & \mbox{on} \;\; \tilde{U} \\ 
            0 & \mbox{on} \;\; \ol{O_4} \\
            \alpha' & \mbox{elsewhere} \;.
            \end{array}\right.
$$
The form $\alpha''$ is smooth and $d_\G \alpha'' = \beta'$. Thus 
$\beta = d_\G (\alpha'' + \alpha)$.\\

\noindent
{\em End of the proof of \lref{quatuor}.}\\

\noindent
This implies that $H^2(\F|_W) = 0$, hence that $H^2(\F) = 0$.
\cqfd 

\begin{prop} The foliated manifold $(M,\F)$ does not admit any \lsste.
\end{prop}

\Pf Suppose on the contrary that $\omega$ is a leafwise \s structure
on $(M,\F)$. As a consequence of \lref{Botticelli}, $\omega = d_\F
\alpha$ for some \td $1$-form $\alpha$. Let $T$ denote the
hypersurface $\partial U_1 \times (-1,1)$. It separates $M$ into two
domains $D^1$ and $D^2$, with $D^1$ contained in $U_1\times
(-1,1)$. Let $p : M \to (-1,1)$ be the restriction of the map $p'$ to 
$M$. For $i=1,2$ and for $t$ in $(-1,1)$, let $D^i_t = D^i\cap p^{-1}(t)$,
let $T_t= T\cap p^{-1}(t)$, let $\omega_t = \omega|_{p^{-1}(t)}$ and
let $\alpha_t = \alpha|_{p^{-1}(t)}$. For every negative~$t$, $D^2_t$
is relatively compact and   
\begin{equation}\label{Monet}
0 < \int_{D^2_t} \omega_t \wedge \omega_t =
\int_{(T_t,\circlearrowright)} \alpha_t\wedge\omega_t \;.
\end{equation}
The orientation implicitly assigned to $T_t$ is that induced from the
volume form $\omega_t \wedge \omega_t$ when $T_t$ is
thought of as the boundary of $D^2_t$. We denote this
orientation by $\circlearrowright$. Similarly, for every positive~$t$
we have
\begin{equation}\label{Cezanne}
0< \int_{D^1_t}\omega_t \wedge \omega_t = \int_{(T_t,
\circlearrowleft)} \alpha_t \wedge \omega_t = -\int_{(T_t,
\circlearrowright)} \alpha_t \wedge \omega_t \;,
\end{equation}
where $\circlearrowleft$ denotes the opposite orientation on $T_t$. Besides, 
the map  
\begin{equation}\label{Degas}
t\mapsto\int_{(T_t,\circlearrowright)}\alpha_t\wedge\omega_t
\end{equation}
is continuous. (Notice that until now we have only used the fact that each
leaf of $\F$ had a trivial second de~Rham cohomology. To be able to say that
the function (\ref{Degas}) is continuous, we need the stronger fact 
that $H^2(\F)=0$.) It follows thus from (\ref{Monet}) and (\ref{Cezanne}) that
\begin{equation}\label{Rubens}
\int_{T_0}\alpha_0\wedge\omega_0 = 0\, .
\end{equation}

Now choose a closed ball $U_1'$ embedded in $F$, containing $p_1$ and
$q_1$, and such that $U'_1 \subset \dim{int} (U_1)$. The same argument
applies to $U_1'$ and shows that   
\begin{equation}\label{Rembrandt}
\int_{T'_0}\alpha_0\wedge\omega_0 = 0\, ,
\end{equation}
for $T'_0 =\partial U'_1\times \{0\}$. Observe that the hypersurface
$T_0 \cup T_0'$ bounds in $F \times \{0\}$ a relatively compact domain
$D_0 = (\dim{int} (U_1) - U'_1) \times \{0\}$. Moreover, we have
$$
0 < \int_{D_0} \omega_0 \wedge \omega_0 = \int_{T_0 \cup T_0'} 
\alpha_0\wedge\omega_0 \;.
$$ 
But (\ref{Rubens}) and (\ref{Rembrandt}) imply that the last expression 
vanishes, a contradiction.
\cqfd

\subsection{Fourth variation}

Let $M$ be an $(n+1)$-dimensional manifold endowed with a codimension one
foliation $\F$ and a Riemannian metric $g$.

\begin{prop}\label{var4} Suppose that the \fm $(M,\F)$ satisfies the 
following condition. For every positive function $\eps : M \to (0,\infty)$,
there exists a nonsymplectic $n$-dimensional submanifold $S$ of $M$ 
whose tangent space $TS$ is $\eps$-close to $T\F$ (by which term is 
meant that $\|X-Y\| < \eps(s)$ for every unit vectors $X$ and $Y$ in 
$T_sS$ and $T_s\F$ respectively and for every $s$ in $S$). Then $(M,\F)$ 
does not support any \lsst that admits a closed extension. 
\end{prop}

\Pf Suppose, on the contrary, that $\Omega$ is a closed $2$-form on
$M$ whose restriction to $T\F\times T\F$ is a \lsste. Consider the
$1$-dimensional distribution $N_\Omega\F = \{X\in TM\,;\,\Omega(X, 
\cdot)=0\}$. The  closed form $\Omega$ restricts to a nondegenerate 
$2$-form, hence a symplectic form, on any hypersurface $S$ transverse 
to $N_\Omega\F$. Now let 
$$
\varepsilon : M \to (0,\infty) : x \mapsto d(T_x\F,(N_\Omega\F)_x)\;,
$$
where $d(T_x\F,(N_\Omega\F)_x) = \dim{inf} \{ \|X - Y\| ; X \in
T_x\F, Y \in (N_\Omega\F)_x, \; \dim{and}\; \|X\|=\|Y\|=1 \}$.
Any submanifold $S$ whose tangent space is $\varepsilon$-close to
$T\F$ is necessarily transverse to $N_\Omega\F$, and therefore
inherits a \sste, contradicting the hypothesis.  
\cqfd

\begin{rmk}{\rm Although the four different arguments for nonexistence
of \lsste s that have been presented here seem independent, it can be 
verified that the last one applies to all our examples. This does not 
remain true if we modify the examples so as to raise their codimension, 
which can be done without affecting existence of \lasste s nor 
nonexistence of \lsste s. 
}\end{rmk}

\bibliographystyle{plain}
\bibliography{mybib}

\end{document}